\documentclass[12pt]{article}
\usepackage{amsmath}
\usepackage{amsfonts}
\usepackage{amssymb}
\usepackage{theorem}

\title{Locally compact groups which are separably categorical structures }
\author{A.Ivanov 
\thanks{
The research is supported by Polish National Science Centre grant DEC2011/01/B/ST1/01406
} 
}

\newtheorem{thm}{Theorem}
\newtheorem{lem}[thm]{Lemma} 
\newtheorem{definicja}[thm]{Definition}
\newtheorem{cor}[thm]{Corollary} 
\newtheorem{prop}[thm]{Proposition} 
\newtheorem{rmrk}[thm]{Remark}

\begin{document} 
\topmargin = 12pt
\textheight = 630pt 
\footskip = 39pt 

\maketitle 

\begin{abstract} 
We describe locally compact groups which are 
separably categorical metric structures.  
\end{abstract}

\section{Introduction}

In this paper we describe separable locally compact groups 
which can be presented as continuous structures 
with  a separably categorical continuous theory. 
The latter means that the group 
is determined uniquely (up to metric isomorphism) by 
its continuous theory among all separable groups. 

Locally compact groups will be considered  
as one-sorted continuous structures 
in the continuous signature 
$$
L=\{ d, \cdot , ^{-1} \},  
$$ 
where $d$ denotes the metric.  
Let us recall that a {\em metric $L$-structure} 
is a complete metric space $(M,d)$ with $d$ bounded by 1, 
where $\cdot$ and $^{-1}$ are uniformly continuous 
operations on $M$ \cite{BYBHU}.    
It is assumed that to the symbols $^{-1}$ 
and $\cdot$ continuity moduli $\gamma_{1}$ 
and $\gamma_2$ are assigned so that 
$$
d(x_1 ,x_2 ) <\gamma_1 (\varepsilon ) \mbox{ implies } d(x^{-1}_1 ,x^{-1}_2 ) <\varepsilon  \mbox{ and } 
$$ 
$$ 
d(x_1 ,x_2 )<\gamma_2 (\varepsilon)  \mbox{ implies } 
max (d (x_1 \cdot y , x_2 \cdot y) , d( y\cdot x_1 , y\cdot x_2 )) < \varepsilon  
\mbox{ for all } y\in M.   
$$  
Note that the latter condition is equivalent to existence of 
$\gamma_3 (\varepsilon )$ (in fact $\gamma_2 (\gamma_2 (\varepsilon ))$) 
such that 
$$ 
d(x_1 ,x_2 )<\gamma_3 (\varepsilon)  \mbox{ implies } 
d( y\cdot x_1 \cdot z, y\cdot x_2 \cdot z )) < \varepsilon  
\mbox{ for all } y,z \in M.   
$$  
We will also assume that the metric is left-invariant. 
We will see in Lemma \ref{1} and Proposition \ref{2} 
below that this assumption is natural. 

We now state the main result of the paper. 
All logic notions appeared in the formulation 
and a proof of the theorem will be given in Section 2.  
We try to make it available for mathematicians 
outside model theory. 

\begin{thm} \label{catlc}
Let $G$ be a locally compact group 
with a left-invariand metric $d \le 1$ 
so that $(G,d)$ is a continuous metric structure.   
Then $(G,d)$ is a separably categorical metric structure  
if and only if there is a compact clopen subgroup $H<G$ 
which is invariant with respect to all metric automorphisms of $G$, 
and the induced action of $Aut(G,d)$ on the coset space $G/H$ 
is oligomorphic. 

In this case if the connected component of the identity 
$G_0$ is a neighbourhood of the identity, 
$H$ can be taken to be $G_0$. 
When this happens or when $d$ is bi-invariant, 
the subgroup $H$ is normal and $G/H$ is an $\omega$-categorical 
discrete group. 
\end{thm} 

We just remind the reader that an action of $G$ on 
a set $\Omega$ is called {\em oligomorphic} 
if for every $n$ the group $G$ has finitely 
many orbits on $\Omega^n$, see Section 2.1 in \cite{ca}.  
Concerning the second part of the statement note that 
in any Lie group the connected component $G_0$  is a neighbourhood of 
the identity (Section 3.0 in \cite{MZ}). 
On the other hand to get an opposite example of 
a separably categorical group just take the direct sum 
of a countably categorical group  with a compact $H$  
such that $H_0$ is not open and $H_0 \not=1$. 
According to Corollary \ref{corlc} below such a group 
has an appropriate metric. 

In the rest of the introduction we discuss 
how powerful the continuous logic approach is 
in the class of locally compact groups. 
We start with  the following lemma 
which shows that continuous logic  
can be applied only to metric groups which 
have bi-invariant metrics. 
This lemma appears as a part of  
Proposition 3.13 in \cite{BY}.

\begin{lem} \label{1} 
Let a group $(G,d)$ be a metric $L$-structure with respect to 
continuity moduli $\gamma_1$ and $\gamma_2$ as above. 
Then $G$ admits a complete bi-invariant metric $d^*$ which 
defines the same topology as $d$.  
\end{lem} 

{\em Proof.} 
We assume that $(G,d)$ is not discrete. 
Let $d^* (x,y) = sup_{u,v} d(u\cdot x\cdot v, u\cdot y\cdot v)$. 
Then  $d^* (x,y)$ is a bi-invariant metric 
with $d(x,y) \le d^* (x,y)$.    
Since for every $\varepsilon$ we have 
$$ 
d(x,y) <\gamma_2 (\gamma_2 (\varepsilon)) \Rightarrow d^* (x,y) <\varepsilon , 
$$ 
each open $d$-ball contains an open $d^*$-ball and vice versa. 
$\Box$ 

\bigskip 

{\em Thus not all locally compact groups can be viewed as metric structures! } 
In fact the continuous logic  approach works only for {\bf SIN}-groups. 

\begin{definicja} (\cite{HMS}, Section 2) 
A topological group $G$ is called a {\bf SIN}-group 
if any neighbourhood of the identity of $G$ 
contains a neighbourhood of the identity  
which is invariant under all inner automorphisms. 
\end{definicja} 

\begin{prop} \label{2} 
A topological group $G$ admitting a complete metric 
can be presented as a continuous metric $L$-structure 
$(G,d)$ under $d$ defining the topology of $G$ if and only if 
$G$ is a {\bf SIN}-group. 
\end{prop} 

{\em Proof.} 
V.Klee has  shown in Section 1 of \cite{klee} that  a 
Hausdorff group $G$ admits a bi-invariant metric if and only if 
it admits a countable complete system of neighbourhoods of the identity 
which are invariant under all inner automorphisms. 
Thus a metrizable group $G$ has a bi-invariant metric if 
and only if $G$ is {\bf SIN}. 
In this case we may assume that 
this metric $d$ is bounded by 1
(otherwise it can be replaced by $\frac{d(x,y)}{d(x,y)+1}$). 
In Section 2 of \cite{klee} it is shown that 
if a completely metrizable topological 
group $G$ has a bi-invariant metric 
then this metric is complete. 
$\Box$

\bigskip

It is worth noting that a locally compact {\bf SIN}-group 
is unimodular. 
Moreover {\bf SIN}-groups can be characterized 
among unimodular locally compact groups as 
those groups $G$ for which the von Neumann algebra 
$VN(G)$ generated by the left regular representation is finite 
(i.e. admits no non-unitary isometry),  see \cite{D}, 13.10.5. 
Thus the class of 
locally compact bi-invariant metric groups 
still contains some interesting families  
(together with obvious examples: locally compact abelian 
groups and discrete groups). 
Moreover it is also known that Polish groups which are 
embeddable into the unitary  
group of a separable finite von Neumann algebra 
(i.e. the algebra has a faithful representation on a separable 
Hilbert space),  admit complete bi-invariant metrics \cite{popa}, 
Section 6.5 (also see Lemma 2.10 of \cite{AM}).

The following corollary of Theorem \ref{catlc} and Lemma \ref{1} 
gives a complete description of locally compact groups 
which can be presented as separably categorical metric structures.  
It substantially restricts the variety of examples mentioned above  
to ones which resemble semidirect products of countably categorical 
discrete groups with compact ones.

\begin{cor} \label{corlc} 
A locally compact group $G$ can be presented as 
a separably categorical metric structure with respect to some metric 
if and only if there is a normal compact clopen subgroup $H<G$ with 
a bi-invariant metric $d$ 
so that $d$ is conjugacy invariant in $G$ 
and the group of automorphisms of $G/H$ 
which are induced by automorphisms of $G$ 
preserving $(H,d)$ is oligomorphic 
(in particular $G/H$ is a countably categorical discrete group). 
\end{cor} 

The proof is given in Section 3. 
In Section 3 we also prove that when $G$ is 
separably categorical and $H$ 
is as in this statement,  stability of the structure induced 
on $G/H$  is equivalent to stability of $G$  with respect to 
continuous formulas of some special type. 

We finally note that although by Lemma \ref{1} 
the metric $d$ in the formulation of Theorem \ref{catlc} 
can be chosen bi-invariant, we do not assume this.  
The choice of $H$ may depend on the metric 
(see Remark \ref{rem}). 

The author is grateful to the referee for helpful remarks.

\section{Preliminaries and the proof} 

\subsection{Necessary preliminaries} 

For convenience of the reader we recall some 
basic definitions from \cite{BYBHU} and \cite{BYU}.  
We keep the signature $L$ together with functions 
$\gamma_1$ and $\gamma_2$  as in Introduction.  
Then atomic formulas are the expressions of the form $d(t_1 ,t_2 )$,
 where $t_i$ are terms (built from functional $L$-symbols). 
They take values from $[0,1]$.   
{\em Statements} concerning metric structures are formulated in the form 
$$
\phi = 0 
$$ 
(called an $L$-{\em condition}), where $\phi$ is a {\em formula}, 
i.e. an expression built from numbers of $[0,1]\cap \mathbb{Q}$ 
and atomic formulas by applications of the following functions: 
$$ 
x/2  \mbox{ , } x\dot- y= max (x-y,0) \mbox{ , } min(x ,y )  \mbox{ , } max(x ,y )
\mbox{ , } |x-y| \mbox{ , } \neg (x) =1-x  \mbox{ , } 
$$ 
$$ 
x\dot+ y= min(x+y, 1) \mbox{ , } x \cdot y \mbox{ , } \sqrt{x}  \mbox{ , }  sup_x \mbox{ and } inf_x . 
$$ 
A {\em theory} is a set of $L$-conditions without free variables 
(here $sup_x$ and $inf_x$ play the role of quantifiers). 
Formulas and statements are interpreted in continuous 
$L$-structures in the natural way.    
For simplicity we often replace expressions of the form 
$\phi \dot- \varepsilon = 0$ with rational $\varepsilon$ 
by $\phi \le \varepsilon$. 

Let $M$ be a continuous metric $L$-structure. 
We define the automorphism group $Aut(M)$ of $M$ to 
be the subgroup of $Iso (M,d)$ consisting of all isometries 
preserving the values of atomic formulas. 
It is easy to see that $Aut(M)$ is a closed subgroup with 
respect to the pointwise convergence topology on $Iso (M,d)$. 

For every $c_1 ,...,c_n \in M$ and $A\subseteq M$ 
we define the $n$-type $tp(\bar{c}/A)$ of $\bar{c}$ over $A$ 
as the set of all $\bar{x}$-conditions with parameters from $A$ 
which are satisfied by $\bar{c}$ in $M$.  
Let $S_n (T_A )$ be the set of all $n$-types over $A$ 
of the expansion of the theory $T$ by constants from $A$. 
There are two natural topologies on this set. 
The {\em logic topology} is defined by the basis consisting of 
sets of types of the form $[\phi (\bar{x})<\varepsilon ]$, 
i.e. types containing some $\phi (\bar{x})\le \varepsilon'$ with 
rational $\varepsilon '<\varepsilon$.   
The logic topology is compact. 

The $d$-topology is defined by the metric 
$$
d(p,q)= inf \{ max_{i\le n} d(c_i ,b_i )| \mbox{ there is a model } N 
\mbox{ with } N\models p(\bar{c})\wedge q(\bar{b})\}. 
$$ 
By Propositions 8.7 and 8.8 of \cite{BYBHU} the $d$-topology is finer 
than the logic topology and $(S_n (T_A ),d)$ is a complete space. 

Definability in continuous structures is introduced as follows. 

\begin{definicja} 
Let $M$ be a continuous metric $L$-structure 
and $A\subseteq M$. 
A uniformly continuous map $P:M^n \rightarrow [0,1]$ 
is called a predicate definable in $M$ over $A$
if there is a sequence $(\phi_k (\bar{x}) :k\ge 1 )$ of 
$L(A)$-formulas such that the maps interpreting 
$\phi_k (\bar{x})$ in $M$ converge to $P(\bar{x})$ 
uniformly in $M^n$. 
\end{definicja} 

It is clear that a definable predicate defines a function 
on $S_n (T_A )$. 
By Propositions 3.4 and 3.10 of \cite{BYU} functions 
$S_n (T_A ) \rightarrow [0,1]$ which are continuous with respect to 
logic topology are precisely those given 
by definable predicates.

A tuple $\bar{a}$ from $M^n$ is {\em algebraic} in $M$ over 
$A$ if there is a compact subset $C\subseteq M^n$ such that 
$\bar{a}\in C$ and the distance predicate $dist(\bar{x},C)$ 
is definable in $M$ over $A$. 
Let $acl(A)$ be the set of all elements algebraic over $A$. 
In continuous logic the concept of algebraicity is 
parallel to that in traditional model theory 
(see Section 10 of \cite{BYBHU}).

\bigskip 

A theory $T$ is {\em separably categorical} if any 
two separable models of $T$ are isomorphic. 
By Theorem 12.10 of \cite{BYBHU} a complete theory $T$ is separably 
categorical if and only if for each $n>0$, every $n$-type $p$ is principal. 
The latter means that for every model $M\models T$, the predicate 
$dist(\bar{x},p(M))$ is definable over $\emptyset$. 

Another property equivalent to separable categoricity states that 
for each $n>0$, the metric space $(S_n (T),d)$ is compact.  
In particular for every $n$ and every $\varepsilon$ there is 
a finite family of principal $n$-types $p_1 ,...,p_m$ so that 
their $\varepsilon$-neighbourhoods cover $S_n(T)$. 

In first order logic a countable structure $M$ is 
$\omega$-categorical if and only if $Aut(M)$ is 
an {\em oligomorphic} permutation group, i.e. 
for every $n$, $Aut(M)$ has finitely many orbits 
on $M^n$. 
In continuous logic we have the following modification.  

\begin{definicja} 
An isometric action of a group $G$ on a metric space $({\bf X},d)$ 
is said to be approximately oligomorphic if for every $n\ge 1$ and $\varepsilon >0$ 
there is a finite set $F\subset {\bf X}^n$ such that 
$$ 
G\cdot F = \{ g\bar{x} : g\in G \mbox{ and } \bar{x}\in F\}
$$
is $\varepsilon$-dense in $({\bf X}^n,d)$. 
\end{definicja} 

\begin{thm} 
(C. Ward Henson, see Theorem 12.10 of \cite{BYBHU} and Theorem 4.25 in \cite{scho}) 
Let $M$ be separable continuous metric structure. 
The following conditions are equivalent. 

(1) The theory of $M$ is separably categorical; 

(2) The $d$-topology on $S_n$ coincides with the logic topology; 

(3) The action of $Aut(M)$ on $(M,d)$ is approximately oligomorphic. 
\end{thm}

\subsection{Proof of Theorem \ref{catlc}}

The sufficiency follows from the description of 
separably categorical structures by approximate 
oligomorphicity of the group of their metric automorphisms.  
Indeed, given $n\ge 1$ and $\varepsilon >0$  
choose a finite $\gamma_2 (\varepsilon )$-net $F_1$ in 
the compact set of all $n$-tuples from  $H$. 
Since $Aut(G,d)$ is oligomorphic on $G/H$ 
there is a finite set $F_2$ of $n$-tuples from 
$G$ which represent all $Aut(G,d)$-orbits 
on the set of all $n$-tuples from $G/H$.  
Then $Aut(G,d) \cdot (F_2 \cdot F_1)$ 
is $\varepsilon$-dense in the set of all 
$n$-tuples of $G$.  

Let us prove the necessity of the theorem. 
We start with the following preliminaries. 
We may assume that $G$ is not discrete. 
There is a non-zero rational number $\rho <1$ such that the $\rho$-ball 
of the unity $B_{\rho} (e)=\{ x\in G: d(x,e)\le\rho\}$ is compact. 
Let $G_{\rho}$ be the subgroup generated by $B_{\rho}(e)$. 
Note that for any $g\in G_{\rho}$ the open ball 
$$
B_{<\rho}( g)=\{ x\in G: d(x,g)<\rho\}= \{ x\in G: d(g^{-1} x,e) <\rho \}
$$
is a subset of $G_{\rho}$; thus $G_{\rho}$ is an open (in fact clopen) subgroup. 
If the connected component of the identity $G_0$ is 
a neighbourhood of the identity, 
we choose $\rho$ such that $B_{\rho}(e) \subseteq G_0$. 
Since $G_0$ is a subgroup of $G$, 
$G_{\rho} \le G_0$.  
Since $G_{\rho}$ is open, $G_{\rho}=G_0$. 
Note that when $d$ is a bi-invariant metric, $G_{\rho}$ 
is a normal subgroup of $G$. 
 
\begin{lem} \label{l8} 
Assume that $G$ is a separably 
categorical  locally compact group. 
Then under the circumstances above 
the predicate $P(x)=d(x,G_{\rho})$ is definable in $G$. 
\end{lem} 

{\em Proof.} 
Since $G_{\rho}$ is preserved under metric automorphisms of 
$G$ the predicate $P(x)$ naturally defines a map 
from the set of $1$-types $S_1$ to $[0,1]$. 
Since $P(x)$ is continuous with respect to   
the $d$-topology on $S_1$, by separable categoricity 
it also continuous with respect to the logic topology. 
The latter exactly means that 
$P(x)$ is definable in $G$. 
$\Box$ 

\bigskip

{\bf Remark.} 
Since $G_{\rho}$ is closed, this lemma states that 
$G_{\rho}$ is a definable set in the sense of continuous logic, 
see Definition 9.16 of \cite{BYBHU}.

Let us also note that since the condition $d(x,e)\le \rho$ defines 
a totally bounded, complete subset in any elementary extension of $G$,   
the set $B_{\rho}(e)$ is a subset of $acl(\emptyset )$ in $G$. 
Thus any 
$B^n_{\rho}(e) = B_{\rho}(e) \cdot B_{\rho}(e)\cdot ....\cdot B_{\rho}(e)$ 
also is a subset of $acl(\emptyset )$. 
In particular $G_{\rho} \subset acl(\emptyset )$.

\begin{lem} 
Under the circumstances above  
there is a natural number $n$ so that 
$G_{\rho}=B^n_{\rho}(e)$. 
In particular $G_{\rho}$ is compact. 
\end{lem} 

{\em Proof.} 
Assume $G_{\rho}\not=B^n_{\rho}(e)$ for all $n\in\omega$.  
Notice that this implies that for every $n$  the $\rho$-neighbourhood 
of $B^n_{\rho}(e)$ does not cover $G_{\rho}$.    
Indeed if $g_1 \in B^n_{\rho} (e)$, $g_2 \in G_{\rho}$ 
and $d(g_1 ,g_2 )<\rho$, 
then $g^{-1}_1 g_2 \in B_{\rho}$ 
and $g_2 = g_1  g^{-1}_1 g_2 \in B^{n+1}_{\rho}$ 
(we use the assumption that $d$ is left-invariant). 
Thus if the $\rho$-neighbourhood 
of $B^n_{\rho}(e)$ covers $G_{\rho}$, 
then we have a contradiction with 
$G_{\rho}\not=B^{n+1}_{\rho}(e)$. 
   
We see that the assumption $G_{\rho}\not=B^n_{\rho}(e)$ 
for all $n\in\omega$ implies that all statements 
$$
sup_{x_1 ...x_n}(min(\rho \dot{-} d(x,x_1 \cdot ....\cdot x_n ), \rho\dot{-}d(e,x_1 ),..., \rho\dot{-}d(e,x_n )))=0
$$ 
are finitely consistent together with $P(x)=0$. 
By Lemma \ref{l8} and compactness of continuous logic 
we obtain a contradiction. 
$\Box$  

\bigskip 

Since $G_{\rho}$ is a characteristic subgroup of $G$ with respect to 
the automorphism group of the metric structure $G$, 
we see that $Aut(G,d)$ acts correctly on $G/G_{\rho}$ 
by permutations of $G/G_{\rho}$. 
Note that $G/G_{\rho}$ is a discrete space with respect to 
the topology induced by the topology of $G$. 

\begin{lem} 
The action of $Aut(G,d)$ on $G/G_{\rho}$ is oligomorphic. 
\end{lem} 

{\em Proof.} 
Since $(G,d)$ is separably categorical, $Aut(G,d)$ 
is approximately oligomorphic on $(G,d)$. 
Thus for every $n$ there is a finite set $F$ of 
$n$-tuples from $G$ such that the set of orbits 
meeting $F$ is $\rho$-dense in $(G,d)$. 
In particular for any $g_1 ,...,g_n \in G$ 
there is a tuple $(h_1 ,...,h_n )\in F$ and 
an automorphism $\alpha \in Aut(G,d)$ 
such that $g^{-1}_i \alpha (h_i )\in G_{\rho}$ 
for all $i\le n$. 
$\Box$ 

\bigskip 

To see that Theorem \ref{catlc} follows from lemmas above 
just take $H$ to be $G_{\rho}$. 

\begin{rmrk} \label{rem} 
{\em We now illustrate that the choice of $H$ may depend on the metric. 
Let $G$ be 
$$
S^1 \times \mathbb{Z}/2\mathbb{Z} \times (\mathbb{Z}/2\mathbb{Z})^{\omega} , 
$$ 
where $S^1 = \{ z\in \mathbb{C}: \parallel z \parallel =1\}$ 
is the circle group and $\parallel x \parallel$ is the Euclidean norm in $\mathbb{C}$. 
Let us consider $G$ with respect to the metric $d_1$ defined as follows:  
$$ 
\mbox{ if } x_1 \mbox{ and } x_2 \mbox{ represent distinct } S^1\mbox{-cosets, then } 
d_1 (x_1  ,x_2 ) = 1 ,  
$$ 
$$ 
\mbox{ if } x_1  \mbox{ and } x_2  \mbox{ represent the same coset, then } 
d_1 (x_1  ,x_2  ) =\frac{1}{4} \parallel 1 - x^{-1}_1 x_2 \parallel .  
$$ 
Applying the proof of Theorem 1 for $\rho = 1/2$ we obtain   
$H_1 = G_{\rho} = B_{\rho}(e) = S^1$. 
Let us correct $d_1$ on the subgroup $S^1 \times \mathbb{Z}/2\mathbb{Z}$ 
as follows: 
$$ 
\mbox{ if } x_1 , x_2 \in S^1 \times \mathbb{Z}/2\mathbb{Z} 
\mbox{ represent distinct } S^1\mbox{-cosets, then } 
d_2 (x_1  ,x_2 ) = \frac{1}{2},  
$$ 
$$ 
\mbox{ if } x_1  \mbox{ and } x_2  \mbox{ represent the same coset, then } 
d_2 (x_1  ,x_2  ) =\frac{1}{8} \parallel 1 - x^{-1}_1 x_2 \parallel .  
$$ 
We obtain a bi-invariant metric $d_2 \le 1$. 
Applying the proof of Theorem 1 for $(G, d_2 )$ and $\rho = 1/2$ we see that 
the corresponding subgroup $H_2$ is different: 
$G_{\rho} = B_{\rho}(e) = S^1 \times \mathbb{Z}/2\mathbb{Z}$. 
}
\end{rmrk}

\section{Locally compact groups, separable categoricity  and stability}

We start with the proof of Corollary \ref{corlc}: 

\begin{quote} 
{\em A locally compact group $G$ can be presented as 
a separably categorical metric structure with respect to some metric 
if and only if there is a normal compact clopen subgroup $H<G$ with 
a bi-invariant metric $d$ 
so that $d$ is conjugacy invariant in $G$ and the group of automorphisms 
of $G/H$ which are induced by automorphisms of $G$ 
preserving $(H,d)$ is oligomorphic 
(in particular $G/H$ is a countably categorical discrete group). }
\end{quote} 

In the proof of sufficency we construct some canonical 
metric $d^*$ corresponding to the topology of $G$, so 
that $(G, d^* )$ is separably categorical. 
In the final part of the section we consider 
the question when $(G,d^* )$ is stable.

\bigskip 

{\em Proof of Corollary \ref{corlc}.} 
If the locally compact group $G$ can be presented as 
a separably categorical metric structure with respect to some metric 
$\le 1$ then by Lemma \ref{1} (and its proof) we may choose an 
equivalent definable metric  $d$ which is bi-invariant. 
Since the automorphism group of $(G,d)$ is still 
approximately oligomorphic, the structure 
$(G,d)$ is separably categorical. 
By Theorem \ref{catlc} we find 
a normal compact clopen subgroup $H<G$ 
so that $G/H$ is an $\omega$-categorical 
discrete group. 
 
To see the sufficiency of the theorem 
assume that $H<G$ is a compact, clopen
normal subgroup. 
Choose a bi-invariant 
metric $d$ as in the formulation. 
We may assume that the $d$-diameter is $\frac{1}{2}$. 
 
We consider the countable $\omega$-categorical group 
$G/H$ with respect to the $\{ 0,1\}$-metric, say $d_0$. 
Then consider $G$ with respect to 
so called {\em wreath product} of metrics 
$d^* = d_0 \mbox{ wrt }d$ \cite{oliynyk}:  
$$ 
\mbox{ if } x_1 \mbox{ and } x_2 \mbox{ represent distinct } H\mbox{-cosets, then } d^* (x_1  ,x_2 ) = 1 ,  
$$ 
$$ 
\mbox{ if } x_1  \mbox{ and } x_2  \mbox{ represent the same coset, then } d^* (x_1  ,x_2  ) = d(1, x^{-1}_1 x_2) . 
$$ 
It is easy to see that $d^*$ is a metric. 
To see that it is a left-invariant metric 
assume that   $x_1$ and $x_2$ represent the same coset with respect to $H$, 
and let 
$$ 
x \cdot x_i = x'_i \mbox{ , } i\in \{ 1,2\} . 
$$ 
Then $x^{-1}_1 x_2 = (x'_1 )^{-1} x'_2$ and  
$d(x_1 ,x_2 ) = d(x'_1 ,x'_2 )$. 
The case when $x_i$ do not represent the same coset is obvious. 
To see that $d^*$ is right-invariant apply the following argument. 
 
We may suppose that   $x_1$ and $x_2$ represent the same coset 
with respect to $H$. 
Let 
$$ 
x_i  \cdot x  = x'_i \mbox{ , } i\in \{ 1,2\} . 
$$ 
Then $x^{-1}x^{-1}_1 x_2 x= (x'_1 )^{-1} x'_2$,  i.e.
$d(e, x^{-1}x^{-1}_1 x_2 x)= d(e, (x'_1 )^{-1} x'_2)$ 
and since $d$ is conjugacy invariant we have 
$d(e, x^{-1}_1 x_2 )= d(e, (x'_1 )^{-1} x'_2)$.  
In particular we see  
$d(x_1 ,x_2 ) = d(x'_1 ,x'_2 )$. 

To finish the proof note that since 
$$ 
\phi (x_1 )^{-1} \phi (x_2 ) = \phi (x^{-1}_1 x_2 ) \mbox{ , where } \phi \in Aut (G), 
$$ 
any automorphism of $G$ which acts on $H$ as an automorphism of 
the continuous structure $(H,d)$, also preserves the metric $d^*$. 
Now the statement that $(G, d^* )$ is 
a separably categorical structure 
follows from sufficiency of Theorem \ref{catlc}.   
$\Box$ 
\bigskip

Given $G$ as in the corollary 
{\em is the group of automorphisms 
of $G/H$ which are induced by automorphisms of $G$ 
preserving $(H,d)$, closed in $Aut(G/H)$? } 
By the final part of the proof above 
such automorphisms always extend to 
automorphisms of the continuous structure $(G,d^*)$.  
Having in mind this question let us consider 
$Aut(G,d^* )$ as a metric group. 
We remind the reader that when $({\bf Y},d)$ is a Polish space  
the corresponding isometry group 
$Iso ({\bf Y})$ is a Polish group 
with respect to the pointwise convergence topology. 
A compatible left-invariant metric can be obtained as follows: 
fix a countable dense set $S=\{ s_i : i\in \{ 1,2,...\} \}$ 
and then define for two isometries $\alpha$ 
and $\beta$ of ${\bf Y}$ 
$$ 
\rho_{S} (\alpha ,\beta )= \sum_{i=1}^{\infty} 2^{-i} min(1, d(\alpha (s_i ),\beta (s_i ))) .
$$ 
Let us fix such a metric on $Aut(G,d^* )$. 

\begin{lem} \label{G/H} 
The group of automorphisms of $G/H$ 
which are induced by automorphisms of $(G,d^*)$ 
is closed in $Aut(G/H)$. 
\end{lem} 

{\em Proof.} 
We denote by $\rho$ the metric $\rho_S$ which we have chosen for $Aut (G,d^* )$. 
Here $S=\{ s_i : i\in \{ 1,2,...\}\}$ is the corresponding countable dense subset of $G$. 
Since $G$ is separable and $H$ is clopen, $G/H$ is countable 
and each $H$-coset contains an element of $S$. 
We order $G/H$ so that the number of any 
$gH$ is determined by the first representative of $gH$ from $S$.   
Note that the metric $\rho_{G/H}$ defined with respect to this 
enumeration of $G/H$ has the following presentation:
$$ 
\rho_{G/H} (\alpha ,\beta )= \sum_{i=1}^{\infty} \{ 2^{-i} :\alpha \mbox{ and } \beta \mbox{ do not agree for the } 
i\mbox{-th element of } G/H \}.
$$ 
Let us consider a sequence $\alpha_i \in Aut (G/H)$, $i\in \omega$, 
converging to some $\alpha \in Aut(G/H)$ with respect to $\rho_{G/H}$. 
We assume that each $\alpha_i$ is induced by some 
$\gamma_i \in Aut(G,d^* )$. 
Let us prove that there is a subsequence of $\{ \gamma_i : i\in \omega \}$ 
converging to some $\gamma \in Aut(G, d^* )$. 
We use the following construction. 

At step $n>0$ we define a partial map $\delta_n$ 
on the $2^{-n}$-net $D_n$ of  the first $n$ cosets of 
the enumeration of $G/H$  so that 
\begin{quote} 
(i) $D_n$ is contained in the minimal initial segment of $S$ 
which covers  the first $n$ cosets of the enumeration of $G/H$
by a $2^{-n}$-net; \\ 
(ii) for any $s\in D_n$ we have $d^* (\delta_n (s), \delta_{n+1} (s))< 2^{-n-1}$ 
(this condition is empty for $n=0$); \\ 
(iii) there is an infinite subsequence of $\alpha_i$ which agree 
with $\delta_n$ on the first $n$ cosets of $G/H$ ;  \\ 
(iv) $\delta_n$ maps $D_n$ into the union of the cosets 
determined by $Im (\alpha_i )$ as above and the subsequence in (iii) 
can be chosen so that for the corresponding $\gamma_i$ 
we have $d^* (\delta_n (s), \gamma_i (s))< 2^{-n-2}$ for all $s\in D_n$. 
\end{quote}  
Let us assume that $D_{n-1}$ and $\delta_{n-1}$ 
are already defined. 
Extend $D_{n-1}$ to some $D_n$ so that (i) is satisfied. 
Take an infinite subsequence of $\alpha_i$ which \\ 
- agree with $\delta_{n-1}$ on the first $n-1$ cosets of $G/H$ and \\ 
- can be extended to a sequence of $\gamma_i$ satisfying 
$$
d^* (\delta_{n-1} (s), \gamma_i (s))< 2^{-n-1} \mbox{ for all }s\in D_{n-1} . 
$$ 
We may also assume that all $\alpha_i$ agree at the $n$-th coset 
of the enumeration of $G/H$.  
Since each $H$-coset is compact, we can find a finite 
$2^{-n-2}$-net, say $U_{n+2}$, in the union of all cosets 
from the image of the first $n$ cosets of $G/H$ 
with respect to our $\alpha_i$-s.   
Any $\gamma_i$ defines the map $D_n \rightarrow U_{n+2}$ 
which takes any $s\in D_n$ to the nearest element 
of $U_{n+2}$. 
Extending $U_{n+2}$ if necessary we may assume 
that these maps are injective. 
Since the set of these maps is finite we find an infinite 
subsequence of $\gamma_i$ defining 
the same map $D_n \rightarrow U_{n+2}$. 
We denote this map by $\delta_n$. 
Then condition (iv) is obvious and condition (ii) for $\delta_n$ 
follows from (iv) for $\delta_{n-1}$ and $\delta_n$  
and the triangle inequality.

Having the sequence $\delta_n$, $n\in \{ 1,2,...\}$, 
we choose $\gamma_n$ so that (iv) is satisfied for 
each pair $\delta_n$ and $\gamma_n$. 
Then for every initial segment of $G/H$ the 
sequence $\gamma_n$ is a Cauchy sequence 
of maps on the union of cosets of this segment 
with respect to the metric $\rho$. 
Thus the sequence $\gamma_n$ converges to 
some $\gamma$ in every union of this form.  
$\Box$ 

\bigskip 
{\em From now on we consider $G/H$ under the structure 
induced by $Aut(G,d^* )$, } which exists by the lemma above. 

We now discuss stability of $(G,d^* )$, where $d^*$ 
is constructed as in the proof of Corollary \ref{corlc}. 
Let us recall the following definition from Section 5.2 of \cite{FHS}.  

A continuous theory $T$ {\bf has the order property} 
if there is a formula $\psi (\bar{x}, \bar{y})$, where 
$\bar{x}$ and $\bar{y}$ are of the same length and sorts,   
and there is a model $M$ of $T$ with $(\bar{a}_i :i \in \omega )\subseteq M$, 
so that 
$$
\psi (\bar{a}_i , \bar{a}_j ) = 0  \Leftrightarrow 
i <j \mbox{ and } 
\psi (\bar{a}_i , \bar{a}_j ) = 1  \Leftrightarrow 
i \ge j . 
$$
By compactness this condition is equivalent 
to the property that for all $n$ and $\delta \in (0,1)$ 
there are $\bar{a}_1 ,..., \bar{a}_n$ such that 
$$
\psi (\bar{a}_i , \bar{a}_j ) \le \delta  \Leftrightarrow 
i <j \mbox{ and } 
\psi (\bar{a}_i , \bar{a}_j ) \ge 1 -\delta  \Leftrightarrow 
i \ge j . 
$$
The theory $T$ is called {\bf stable} if it does not have the order property. 
In fact this coincides with the definition of stability given in 
Section 8 of \cite{BYU}.  

We will concentrate on formulas $\psi$ which behave 
as slow functions. 

\begin{definicja} 
Let $r\in [0,1]$. 
We will say that a continuous formula $\psi (\bar{x})$ 
is $r$-{\bf slow} if there do not exist 
tuples $\bar{b}$ and $\bar{b}'$  
with $max (d(b_i ,b'_i )) \le r$ so that 
$\psi (\bar{b})=0$ and $\psi (\bar{b}')=1$. 
\end{definicja} 

In the situation of $(G,d^* )$, where $d^*$ 
is constructed as in the proof of Corollary \ref{corlc}, if $b$ and $b'$ 
represent the same cosets with respect to $H$ then 
the distance between them is $\le 1/2$.  
In particular if $n$-tuples 
$\bar{b}$ and $\bar{b}'$ 
represent the same cosets with respect to $H$ 
in an appropriate power $G^k$ then 
they are not distant at $> 1/2$ under 
the corresponding $max$-metric defined by $d^*$. 
Thus  in this situation a formula $\psi (\bar{x})$ 
is $\frac{1}{2}$-slow if and only if 
there do not exist tuples $\bar{b}$ and $\bar{b}'$  
which represent the same cosets with respect to $H$ 
so that $\psi (\bar{b})=0$ and $\psi (\bar{b}')=1$. 

\begin{thm} 
Let $G$ be a locally compact group which can be presented 
as a separably categorical metric structure and 
let $H$ and $d^*$ be as in Corollary \ref{corlc} and its proof.   

Then the continuous theory $Th (G,d^* )$ does not have $\frac{1}{2}$-slow 
formulas with the order property if and only if 
the elementary theory of the structure $G/H$ is stable. 
\end{thm} 

{\em Proof.} 
Assume that  $Th (G,d^*)$ has the order property 
witnessed by a $\frac{1}{2}$-slow formula 
$\psi (\bar{x}, \bar{y})$ and an $\omega$-sequence 
$(\bar{a}_i , i\in \omega )$. 
We know that $\psi$ is uniformly continuous. 
Thus replacing $\psi$ by appropriate roots 
$\psi^{2^{-l}}$ if necessary, we can arrange that 
for any tuples $\bar{b}_1, \bar{b}'_1 , \bar{b}_2$ 
and $\bar{b}'_2$ of length $ |\bar{a}_i |$ with 
$max (\psi (\bar{b}_1 ,\bar{b}_2 ), d^* (\bar{b}_1 \bar{b}_2 , \bar{b}'_1 \bar{b}'_2 ))< 1/2$
the value $\psi (\bar{b}'_1 ,\bar{b}'_2 )$ is less than 1. 
In particular defining $\psi (\bar{x}, \bar{y})$ 
on the discrete set $(G/H)^k$  by the minimal values 
on cosets corresponding to $\bar{x}$ and $\bar{y}$ 
we see that the relation $\psi (\bar{x}, \bar{y}) < 1/4$
has the order property on $G/H$ in the classical sense. 
Since it is invariant under the oligomorphic automorphism 
group of the structure on $G/H$ induced by $Aut(G, d^* )$ 
we see that the theory of $G/H$ is not stable.  

The opposite direction is easy. 
Having a relation $\theta (\bar{x}, \bar{y})$ with 
the order property with respect to the structure on 
$G/H$ we may define it on $G$ so that it takes constant 
values on cosets: $0$ or $1$. 
Since the elements of distinct cosets 
are distant by $1$, we see that  
$\theta$ (viewed as a continuous formula) is continuous in $(G,d^* )$. 
It is also clear that it is $Aut(G,d^* )$-invariant. 
Applying the argument of Lemma \ref{l8} we see that $\theta$ is 
definable as a continuous predicate (by a sequence of formulas). 
Take a close approximation of $\theta$, say $\theta'$, which is a formula. 
If $\theta'$ does not have values $0$ or $1$, we can correct it 
as $2\theta' - \varepsilon$ for sufficiently small rational $\varepsilon$.   
As a result we obtain a $\frac{1}{2}$-slow formula with 
the continuous order property.  
$\Box$ 

\bigskip 

It is worth noting that the condition of 
$\frac{1}{2}$-slowness is essential in the theorem. 
For example an $\aleph_0$-categorical extra-special 
$p$-group $G$ has a finite definable subgroup $H$ 
so that $G/H$ is elementary abelian and thus is stable. 
On the other hand it is shown in \cite{F} that extra-special 
$p$-groups are not stable.    
Our theorem shows that if we introduce appropriate 
$d^*$ (which defines the discrete topology in this case) 
the first-order formula defining the order property 
is not $\frac{1}{2}$-slow.

\bigskip

1. Institute of Mathematics, Wroc{\l}aw University, pl.Grunwaldzki 2/4, 50-384 Wroc{\l}aw, Poland, \\ 

2. Institute of Mathematics, Silesian University of Technology, ul.Kaszubska 23, 44-100 Gliwice, Poland \\ 

e-mail: iwanowaleksander@gmail.com  
  

\begin{thebibliography}{99} 
\bibitem{AM} Ando, H., Matsuzawa, Y.: On Polish groups of finite type. 
Publ. Res. Inst. Math. Sci. {\bf 48}, 389 - 408 (2012)
\bibitem{BY} Ben Yaacov, I.: Definability of groups in $\aleph_0$-stable metric structures.  
J.Symbolic Logic, {\bf 75}, 817 - 840 (2010)  
\bibitem{BYBHU} Ben Yaacov, I., Berenstein, A., Henson, W., Usvyatsov, A.:  
{\em Model theory for metric structures}. In:  Model theory with 
Applications to Algebra and Analysis, v.2  (Chatzidakis, Z., Macpherson, H.D.,   
Pillay, A., Wilkie, A. eds.). London Math. Soc. Lecture Notes, v.350, 
pp. 315 - 427, Cambridge University Press (2008)   
\bibitem{BYU} Ben Yaacov, I., Usvyatsov A.: Continuous first order logic and local stability. 
Trans. Amer. Math. Soc. {\bf  362}, 5213 - 5259 (2010) 
\bibitem{ca} Cameron, P.:  Oligomorphic Permutation Groups.  
London Math. Soc. Lecture Notes, vol. 152,  
Cambridge University Press (1990) 
\bibitem{D} Diximier, J.: Les $C^*$-alg\'{e}bres et leurs Repr\'{e}sentations. 
Gauthier-Villas, Paris (1969)
\bibitem{FHS} Farah, I., Hart B., Sherman D.: Model theory of operator algebras II: 
Model theory. Israel J. Math. {\bf 201}, 477 - 505 (2014) 
\bibitem{F} Felgner, U.: On $\aleph_0$-categorical extra-special $p$-groups. 
Comptes Rendus de la Semaine d'\'{E}tude en Th\'{e}orie des Mod\'{e}les.  
Logique et Analyse {\bf 18}, no. 71-72, 407 - 428  (1975)
\bibitem{HMS} Hofmann, K.H., Morris, S., Stroppel, M.: Varieties of topological groups, 
Lie groups and SIN groups. Coll. Math., 70, 151 - 163 (1996) 
\bibitem{klee} Klee, V.L., Jr.: Invariant metrics in groups (solution of a problem of Banach). 
Proc. Amer. Math. Soc. {\bf 3}, 484 - 487(1952)  
\bibitem{MZ} Montgomery, D., Zippin, L.: Topological Transformation Groups. 
Interscience Publishers, New York (1955) 
\bibitem{oliynyk} Oliynyk, B.: Isometry groups of wreath products of metric spaces.  
Algebra and Discr. Math. 4, 123 - 130 (2007)  
\bibitem{popa} Popa, S.: Cocycles and orbit equivalence superrigity 
for malleable actions of $w$-rigid groups. Invent. Math. {\bf 170}, 243 - 295 (2007)  
\bibitem{scho} Schoretsanitis, K.:  Fra\"{i}ss\'{e} Theory for Metric Structures.   
PhD thesis,  University of Illinois at Urbana-Champaign (2007) \\ 
(available at http://www.math.uiuc.edu/~henson/cfo/metricfraisse.pdf)    

\end{thebibliography}
\end{document}